\documentclass[12pt,a4paper,twoside]{article}

\usepackage[leqno]{amsmath}
\usepackage{amssymb}
\usepackage{amsthm}

\usepackage{array}

\usepackage{enumerate}

\def\nc{\newcommand}

\let\phi\varphi
\let\epsilon\varepsilon

\let\subset\subseteq
\let\setminus\smallsetminus

\nc{\fa}{a.e.}
\nc{\eg}{e.g.}
\nc{\ie}{i.e.}
\nc{\wolog}{w.l.o.g.}
\nc{\wrt}{w.r.t.}

\newcommand{\zahlen}{\mathbb}

\newcommand{\N}{{\zahlen N}}
\newcommand{\Q}{{\zahlen Q}}
\newcommand{\R}{{\zahlen R}}
\newcommand{\Z}{{\zahlen Z}}

\newcommand{\on}{\operatorname}

\nc{\graph}{\on{graph}}
\nc\id{\on{id}}

\nc{\ch}[1]{\mathbf{1}_{#1}}
\nc{\restr}{|}

\nc{\trand}{\partial}
\nc{\cl}{\on{cl}}
\nc{\interior}{\on{int}}
\nc{\clos}[1]{\overline{#1}}

\newcommand{\betrag}[1]{\lvert #1 \rvert}
\newcommand{\norm}[1]{{\| #1 \|}}

\nc{\ball}{B}
\nc{\oball}{\ball}
\nc{\cball}{\overline B}
\nc{\rball}{\trand B}

\nc{\dist}{\on{dist}}
\nc{\diam}{\on{diam}}
\nc{\vol}{\on{vol}}
\nc{\area}{A}

\nc{\Leb}{\mathcal{L}}
\nc{\Haus}{\mathcal{H}}
\nc{\hHaus}{\Haus^{n-1}}

\nc{\aplimsup}{\operatorname*{ap\,lim\,sup}}
\nc{\apliminf}{\operatorname*{ap\,lim\,inf}}
\nc{\esssup}{\operatorname*{ess\,sup}}

\nc{\dR}{\mathrm{dR}}

\nc{\Hhmr}{H_{n-1}(M,\R)}
\nc{\Hmmr}{H_m(M,\R)}
\nc{\Hhmz}{H_{n-1}(M,\Z)}
\nc{\Hmmz}{H_m(M,\Z)}
\nc{\Himr}{H_1(M,\R)}
\nc{\Himzr}{H_1(M,\Z)_\R}
\nc{\Hmdrm}{H_{\dR}^m(M)}
\nc{\Hidrm}{H_{\dR}^1(M)}

\nc{\loc}{\mathrm{loc}}

\nc{\normal}{\mathbf N}
\nc{\hnormal}{\normal_{n-1}}
\nc{\lnormal}{\normal^\loc}
\nc{\nmm}{\normal_m(M)}
\nc{\lnmm}{\lnormal_m(M)}
\nc{\nmrn}{\normal_m(\R^n)}
\nc{\hlnormal}{\normal^\loc_{n-1}}
\nc{\cnormal}{\normal^c}
\nc{\flt}{\mathbf F}
\nc{\rflat}{\mathbf F}
\nc{\hflat}{\rflat_{n-1}}
\nc{\fmm}{\rflat_m(M)}
\nc{\lrflat}{\rflat^\loc}
\nc{\integral}{\mathbf I}
\nc{\lintegral}{\integral^\loc}
\nc{\rect}{\mathcal R}
\nc{\lrect}{\rect^\loc}
\nc{\hrect}{\rect_{n-1}}
\nc{\hlrect}{\rect^\loc_{n-1}}
\nc{\lbv}{BV^\loc}

\nc{\rcyc}{\mathbf Z}
\nc{\hrcyc}{\rcyc_{n-1}}
\nc{\hrcycm}{\hrcyc(M)}
\nc{\mrcycm}{\rcyc_m(M)}
\nc{\rbound}{\mathbf B}
\nc{\hrbound}{\rbound_{n-1}}
\nc{\hrboundm}{\hrbound(M)}
\nc{\mrboundm}{\rbound_m(M)}

\newcommand{\rand}{\partial}
\nc{\rrand}{\partial^*}

\newcommand{\mass}{\mathbf{M}}

\newcommand{\ldbr}{\mathopen{[\mspace{-3mu}[}}
\newcommand{\rdbr}{\mathclose{]\mspace{-3mu}]}}
\newcommand{\current}[1]{\ldbr #1 \rdbr}

\newcommand{\inner}{\haken}
\def\haken{\mathbin{\vrule height0.6em \vrule height0.4pt width0.6em}}

\newcommand{\push}[1]{{#1}_{\#}}

\nc{\spt}{\on{spt}}
\nc{\reg}{\on{reg}}
\nc{\sing}{\on{sing}}

\nc{\stab}{\mathcal S}
\nc{\PD}{\mathrm{PD}}

\newcommand{\theoremshape}{\slshape}

\newtheoremstyle{myplain}%
    {\topsep}{\topsep}{\theoremshape}%
    {}{\bfseries}{.}{ }%
    {\thmnumber{#2}\thmname{ #1}\thmnote{ (#3)}}
\newtheoremstyle{mydefinition}%
    {\topsep}{\topsep}{\upshape}%
    {}{\bfseries}{.}{ }%
    {\thmnumber{#2}\thmname{ #1}\thmnote{ (#3)}}
\newtheoremstyle{myremark}%
    {\topsep}{\topsep}{\upshape}%
    {}{\slshape}{{.}}{ }%
    {\thmname{#1}\thmnote{ (#3)}}
\newtheoremstyle{namethm}%
    {\topsep}{\topsep}{\theoremshape}%
    {}{\bfseries}{.}{ }%
    {\thmnumber{#2}\thmnote{ #3}} 
\newtheoremstyle{nonumnamethm}%
    {\topsep}{\topsep}{\theoremshape}%
    {}{\bfseries}{.}{ }%
    {\thmnumber{}\thmnote{#3}}
\newtheoremstyle{nurnummer}%
    {\topsep}{\topsep}{\upshape}%
    {}{\bfseries}{}{ }%
    {\thmnumber{#2}} 
\theoremstyle{myplain}

\newtheorem{theorem}{Theorem}[section]
\newtheorem{lemma}[theorem]{Lemma}
\newtheorem{proposition}[theorem]{Proposition}
\newtheorem{corollary}[theorem]{Corollary}

\theoremstyle{mydefinition}
\newtheorem{definition}[theorem]{Definition}

\theoremstyle{myremark}
\newtheorem*{remark}{Remark}

\theoremstyle{namethm}

\theoremstyle{nurnummer}

\theoremstyle{nonumnamethm}
\newtheorem*{citthm}{}

\nc{\frage}{\marginpar{\Large \sffamily \bfseries ?}}
\nc{\quer}{\bar}
\nc{\dach}{\hat}
\nc\eucl{\mathrm{eucl}}
\def\bref(#1){(\ref{(#1)})}
\nc{\levels}{\mathcal L}
\nc{\leaves}{\mathcal N}
\nc{\hleaves}{\leaves_{n-1}}
\nc{\Kern}[1]{K(#1)}

\begin{document}

\title{Differentiability of the stable norm in codimension one}
\author{Franz Auer%
          \footnote{Research supported by the Communaut\'e fran\c caise de 
                    Belgique, through an Action de Recherche Concert\'ee de 
                    la Direction de la Recherche Scientifique}  
       \ and Victor Bangert}
\date{November 30, 2004}

\maketitle

\begin{abstract}
The real homology of a compact, $n$-dimensional Riemannian manifold $M$
is naturally endowed with the stable norm.
The stable norm of a homology class is the
minimal Riemannian volume of its representatives.
If $M$ is orientable the stable norm on $H_{n-1}(M,\R)$ is a homogenized
version of the Riemannian $(n{-}1)$-volume.
We study the differentiability properties of the stable norm at points
$\alpha \in \Hhmr$. 
They depend on the position of $\alpha$ with respect to the integer lattice
$\Hhmz$ in $\Hhmr$.
In particular, we show that the stable norm is differentiable at $\alpha$ if
$\alpha$ is totally irrational.
\end{abstract}

\section{Introduction}

On every compact Riemannian manifold $M$ the real homology vector spaces 
$\Hmmr$ are endowed with a natural norm,
called \emph{stable} or \emph{mass norm}.
The stable norm $\stab(\alpha)$ of $\alpha \in \Hmmr$
is defined as
the infimum of the Riemannian $m$-volumes of real singular cycles
representing $\alpha$.
Equivalently, $\stab(\alpha)$ can be defined as the minimum of the masses
of closed $m$-currents representing $\alpha$.
The term ``stable norm'' was coined by M.~Gromov, cf.~\cite[Chapter 4]{Gr}.
The concept itself was introduced prior to this by H.~Federer 
\cite[4.15]{Fe2}.
If $m=1$, or if $M$ is $n$-dimensional and orientable and $m=n-1$, then
the stable norm is a homogenized version of the Riemannian length or
$(n{-}1)$-volume functional, respectively.
Here, homogenization is performed with respect to $\Z^b$ acting as group
of deck transformations on the Abelian covering of $M$,
where $b=b_1(M)$ denotes the first Betti number.

We study differentiability properties of the stable norm $\stab$ in the
codimension one case, \ie\ in the case $m=n-1$.
At a point $\alpha \in \Hhmr$ the existence of two-sided 
directional derivatives
of $\stab$ at $\alpha$ depends on the position of $\alpha$ with respect to 
the integer lattice $\Hhmz$ in $\Hhmr$.

\begin{theorem}\label{theorem}
  Let $M$ be a compact, orientable, $n$-dimensional Riemannian manifold 
  and $\stab \colon \Hhmr \to \R_{\ge0}$ the associated stable norm on 
  $\Hhmr$.
  If $\alpha \in \Hhmr$, let $V(\alpha)$ denote the smallest linear
  subspace of $\Hhmr$ that is spanned by integer classes and contains $\alpha$.
  Then 
  the restriction of $\stab$ to $V(\alpha)$ is differentiable at $\alpha$.
\end{theorem}

The extremal cases are that $\alpha$ is rationally independent, in which case
$V(\alpha) = \Hhmr$ and $\stab$ is differentiable at $\alpha$, and the case
that the direction of $\alpha\ne0$ is rational, 
in which case $\dim V(\alpha) = 1$
and the claim of Theorem~\ref{theorem} is obvious.

Due to the convexity and homogeneity of $\stab$ the claim of 
Theorem~\ref{theorem} can be stated in the following alternative form.

{\sl
The tangent cone of the unit ball 
$B := \big\{ \beta \in \Hhmr \bigm| \stab(\beta) \le 1 \big\}$
at $\alpha \in \rand B$ splits as a product, with one factor a
hyperplane in~$V(\alpha)$.}

There is strong evidence that this result is optimal, in the sense that for
a large set of Riemannian metrics on an $n$-torus $T^n$ the stable norm on 
$H_{n-1}(T^n,\R) \simeq \R^n$ is two-sided differentiable precisely in the 
directions covered by Theorem~\ref{theorem}, cf.~\cite{Ba1}, \cite{Se2},
\cite[10.4]{CL}.
On the other hand, for flat metrics on $T^n$ the stable norm on 
$H_{n-1}(T^n,\R)$ is induced by a scalar product.
For this and other explicit examples, see \cite[4.15]{Fe2}.

In the case of the $2$-torus $T^2$ the theorem is proved in \cite{Ba1}, see
also \cite{Au} and \cite{Ma}.
For closed, orientable Riemannian surfaces $F$ of genus $s > 1$ the 
boundary structure of the stable norm ball 
$B \subset H_1(F,\R) \simeq \R^{2s}$ is studied in \cite{Mas}
and \cite{Mas2}.
In particular, in this case  Theorem~\ref{theorem} follows from
\cite[Corollary 3]{Mas2}.
W.\ Senn \cite{Se1}, \cite{Se2}, 
\cite{Se3} proved results analogous to
Theorem~\ref{theorem} for  $\Z^n$-periodic nonparametric variational problems.
The case of the stable norm on $H_1(T^n,\R)$ for $n >2$ is considerably more
subtle, see \cite{BIK}.

Although the basic idea of our proof for Theorem~\ref{theorem} is simple,
we meet some complications that are caused by the lack of regularity of the
objects involved.
Here, we give a rough sketch of the proof that does not attend to such 
subtleties.
The proof is based on the duality between real homology and cohomology
realized by flat cycles and cocycles from Geometric Measure Theory.
Dual to the stable norm on $\Hmmr$ we have the comass norm on the de Rham
cohomology vector space $\Hmdrm$.
Here the comass of a cohomology class $l \in \Hmdrm$ is the infimum of the
maximum norms of smooth closed $m$-forms representing $l$.
An important point in the proof is the existence of a bounded, measurable,
weakly closed $m$-form $\lambda$ that represents~$l$, realizes this 
infimum and that is defined by a process of differentiation from a flat
cocycle representing $l$.
This is due to J.\,H. Wolfe, cf.\ \cite{Wh}.

Differentiability of $\stab\restr V(\alpha)$ at $\alpha \in \Hhmr$ 
means that $l_1(\beta) = l_2(\beta)$ whenever
$\beta \in V(\alpha)$ and
$l_1, l_2 \in H^{n-1}_{\rm dR}(M)$ are subderivatives of $\stab$ at $\alpha$.
Hence, in order to prove Theorem~\ref{theorem}, we will show that
\begin{equation}\label{e:1.1}
(l_1 - l_2)(\beta) = 0 \,.
\end{equation}
We let $\lambda_1, \lambda_2$ denote the $(n{-}1)$-forms representing
$l_1, l_2$ mentioned above.
We choose a smooth closed $1$-form $\eta$ 
that represents the Poincar\'e dual
of $\beta$, so that
$$
 (l_1 - l_2)(\beta) = \int_M \eta \wedge (\lambda_1 - \lambda_2) \,.
$$
Finally, there exists a closed $(n{-}1)$-current $T$ representing $\alpha$
and of minimal mass $\mass(T) = \stab(\alpha)$.
Then the assumption that $l_1$ and $l_2$ are subderivatives 
of $\stab$ at
$\alpha$ can be used to prove that $\lambda_1$ and $\lambda_2$ coincide
at Lebesgue almost all points
of $\spt T$.
So, in order to prove \eqref{e:1.1} it suffices to show that
$$
 \int_{M\setminus \spt T} \eta \wedge (\lambda_1 - \lambda_2) = 0 \,.
$$
According to \cite{AB'} the current $T$ can be represented as a measured
lamination by minimizing hypersurfaces (possibly with singularities).
The connected components of $M \setminus \spt T$ are called the 
gaps of $T$ and it remains to prove that
\begin{equation}\label{e:1.2}
 \int_G \eta \wedge (\lambda_1 - \lambda_2) = 0
\end{equation}
for every gap $G$ of $T$.
Using the fact that $\beta \in V(\alpha)$, one can see that $\eta$ is exact
on $G$, $\eta = dg$ for some function $g \in C^\infty(G,\R)$.
Hence, in the weak sense, we have
\begin{equation}\label{e:1.3}
  \eta \wedge (\lambda_1 - \lambda_2) =  d \big(g\,(\lambda_1 - \lambda_2)\big)
\end{equation}
on $G$.
Now, one would like to  integrate \eqref{e:1.3} over $G$ using
Stokes's theorem and to conclude that the boundary terms vanish since
$\lambda_1 = \lambda_2$ on $\rand G$. This would prove \eqref{e:1.2}.
Actually, one has to be more careful at this point since on the one hand
$\overline G$ is not a compact domain with smooth boundary, 
and on the other
hand $\lambda_1$ and $\lambda_2$ need not be defined on $\rand G$.

In the course of the proof we obtain the following result which may be of 
independent interest.
For the notions in this statement
see the beginnings of
Sections \ref{s:3.1n} and~\ref{s:4n}.

\begin{citthm}[\ref{reg-reg} Theorem]
 Suppose that the flat cocycle $L$ is a calibration and that 
 the closed rectifiable current
 $T \in \rect_{n-1}(M)$ is calibrated by $L$.
 Then the singular set of $T$ is contained in the singular set of $L$.

 In particular, the union of the singular sets of all closed 
 $T \in \rect_{n-1}(M)$ 
 calibrated by $L$ is a Lebesgue null set.
\end{citthm}

As the preceding statement shows we have to use notions and results
from Geometric Measure Theory.
In order to make the article reasonably comprehensive we give
definitions for most of the concepts that are used in an essential way.
These are of a functional analytic nature and easily comprehensible.
In 
Section~\ref{s:2n} we treat the existence of 
mass-minimizing currents in real homology classes.
Section~\ref{s:3n} develops the dual theory of flat cocycles, calibrations
and their representation by weakly closed $L^\infty$-forms.
In Section~\ref{s:4n} we specialize to the codimension one case and
prove Theorem~\ref{reg-reg} above.
In Section~\ref{s:5n} we summarize results on the structure of mass-minimizing
currents in codimension one real homology classes.
These are formulated for the lift of the current to the smallest
(infinite) covering space on which this lift bounds.
Finally, in Section~\ref{s:6n}, the proof of Theorem~\ref{theorem}
is completed.

\section{Currents and the Stable Norm}\label{s:2n}

\subsection{The stable norm}\label{s:2.1n}

Throughout the paper, $M$ will denote an oriented, $n$-dimensional
Riemannian manifold and $m$ an integer, $0 \le m \le n$.
In this subsection we assume in addition that $M$ is compact.

The \emph{stable} (or \emph{mass}) \emph{norm} $\stab(\alpha)$ of
a real homology class $\alpha \in \Hmmr$ is defined as the infimum of the
Riemannian volumes of all real Lipschitz cycles $c = \sum r_i \sigma_i$
representing $\alpha$, see \cite[4.C]{Gr}.
Here the volume $\vol_m(c)$ of~$c$ is $\sum \betrag{r_i} \vol_m(\sigma_i)$,
where $\vol_m(\sigma_i)$ denotes the $m$-dimensional Riemannian volume of the
Lipschitz simplex $\sigma_i \colon \Delta^m \to M$.
To see that $\stab(\alpha)>0$ if $\alpha \ne 0$, note that there exists a 
de Rham cohomology class $\beta \in \Hmdrm$ such that
$$
 0 < [\alpha,\beta] = \int_c \omega  \le \vol_m(c)\, \norm{\omega}_\infty\,,
$$
whenever $c$ represents $\alpha$ and the closed $m$-form $\omega$ represents
$\beta$.

In general, there will not exist a real Lipschitz cycle $c$ in a given
homology class $\alpha \in \Hmmr$ such that the volume of $c$ equals the 
stable norm $\stab(\alpha)$ of $\alpha$.
Geometric Measure Theory provides an appropriate notion of weak solution
to this problem -- the normal currents.

\subsection{Normal and locally normal currents} \label{s:2.2n}

In the following we do not assume that $M$ is compact.
We consider the chain complex
$$
 \rand\colon \big(\Omega_0 ^{m+1} M\big)^* \to \big(\Omega_0 ^{m}M\big)^*
$$
that is dual to the complex
$$
 d \colon \Omega_0 ^{m}M \to \Omega_0 ^{m+1} M
$$
given by the spaces $\Omega_0 ^{m}M$ of smooth $m$-forms $\omega$ with compact
support $\spt \omega$ in $M$ and the exterior derivative.
On $\Omega_0^m M$ we consider the \emph{comass norm}
$$
 \norm \omega_\infty := \max_{x \in M} \norm{\omega_x} \,,
$$
where $\norm{\omega_x}$ denotes the (pointwise) \emph{comass norm} of 
$\omega_x \in \Lambda^m TM_x$,
$$
 \norm{\omega_x} = \max \big\{\omega_x(e_1,\dots,e_m) \bigm|
    e_i \in TM_x \text{ and } \betrag{e_i} \le 1 \text{ for } 1 \le i \le m
      \big\} \,.
$$
The \emph{mass} $\mass(T) \in [0,\infty]$ of $T \in \big(\Omega_0 ^{m}M\big)^*$
is defined by
$$
 \mass(T) = \sup\big\{T(\omega) \bigm| \omega \in \Omega_0^mM 
     \text{ and } \norm\omega_\infty \le 1 \big\} \,.
$$ 
We say that $T$ is of \emph{locally finite mass} if the open sets $U \subset M$
such that
$$
 \mass_U(T) = \sup\{T(\omega) \bigm| \omega \in \Omega_0^m U
         \text{ and }
     \norm\omega_\infty \le 1 \big\} < \infty
$$
cover $M$.
In this case $T$ is representable by integration, i.e., there exist a 
positive
Radon measure $\mu_T$ and a $\mu_T$-measurable unit $m$-vector field
$\vec T$ such that
$$
 T(\omega) = \int \langle \vec T, \omega \rangle \, d\mu_T
$$
for all $\omega \in \Omega_0^mM$.
Clearly $\mu_T(U) = \mass_U(T)$ if $U \subset M$ is open.
If $T$ is of locally finite mass and $f \colon M \to \R$ is 
locally $\mu_T$-integrable
then $T \inner f$, defined by 
$$
 \big(T \inner f\big)(\omega) = 
  \int f\, \langle \vec T, \omega\rangle \, d\mu_T \,,
$$
is of locally finite mass. 
If $A \subset M$ is $\mu_T$-measurable  
one sets
$$
     T \inner A := T \inner \chi_A\,.
$$
{\sloppy
\begin{definition}\label{2.1n}
 A \emph{locally normal current $T$ on $M$} is an element of 
 $\big(\Omega_0^m M\big)^*$ such that both $T$ and $\rand T$ have locally
 finite mass.
 We let $\lnormal_m(M)$ denote the set of locally normal currents on $M$.
 A \emph{normal current} is a locally normal current $T$ whose support
 $\spt T$ ($= \spt \mu_T$) is compact.
 We let $\normal_m(M)$ denote the set of normal currents.
\end{definition}

}
Note that the mass $\mass$ is a norm on the $\R$-vector space $\normal_m(M)$.
Here are the most important classes of examples of (locally) normal currents.

\begin{enumerate}
\item \emph{Lipschitz chains}.
If $c = \sum r_i \sigma_i$ is a real Lipschitz $m$-chain in $M$, then
$\current c \in \normal_m(M)$ is defined by
$$
 \current c (\omega) = \int_c \omega = 
  \sum r_i \int_{\Delta_m} \sigma_i^* \omega \,.
$$
Stokes's Theorem implies $\rand \current c = \current{\rand c}$.
Moreover 
$$
 \mass\big(\current c\big) \le \vol_m(c) \,.
$$
Contrary to what is stated in \cite[4.17]{Gr}, inequality can occur.
If parts of the singular simplices cover each other with different
orientations, these parts add to $\vol_m(c)$, while they cancel in 
$\current c$.
If $m < n$ this situation does not occur for generic chains.

\item \emph{Smooth currents}.
If $\eta \in \Omega^{n-m} M$ is a smooth $(n{-}m)$-form then
$T_\eta = \current M \inner \eta \in \lnormal_m(M)$ is defined by 
$$
 T_\eta(\omega) = \int_M \eta \wedge \omega \,.
$$
Then one has
$$
 \rand T_\eta = (-1)^{n-m+1} T_{d\eta} \,.
$$
If $m \in \{0,1,n-1,n\}$ then
$$
 \mass(T_\eta) = \int_M \norm{\eta(x)} \, d\vol_n(x) \,.
$$
A similar equality is also true for $1 < m < n-1$ if one
uses an appropriate norm on $(n{-}m)$-covectors 
(namely the norm dual to the comass norm on $m$-vectors with respect to the
duality induced by the wedge product and the volume form).
Obviously, one has $T_\eta \in \normal_m(M)$ if and only if 
$\eta \in \Omega_0^{n-m} M$.
\end{enumerate}

\subsection{Flat norm}\label{s:2.3n}

With respect to the mass norm neither the subspace of smooth currents 
with compact support nor the subspace of Lipschitz chains is dense 
in $\normal_m(M)$. 
However, this is true with respect to the flat norm $\flt$,
which is weaker than the mass norm $\mass$.

\begin{definition}\label{2.2n}
 The \emph{flat norm} $\rflat(T)$ of $T \in \normal_m(M)$ is defined by
 $$
  \rflat(T) = 
  \inf \big\{\mass(T+\rand S) + \mass(S) \bigm| S \in \normal_{m+1}(M) \big\}
  \,.
 $$
\end{definition}
One can show that
$$
 \rflat(T) = \sup \big\{T(\omega) \bigm|
     \omega \in \Omega_0^mM,\ \norm\omega_\infty \le 1,\ 
               \norm{d\omega}_\infty \le 1 \big\}\,,
$$
cf.\ \cite[4.1.12]{Fe}, and that the two subspaces mentioned above are
$\rflat$-dense in $\normal_m(M)$, cf.\ \cite[4.1.18 and 4.1.23]{Fe}.

\subsection{Homologically mass-minimizing currents}\label{s:2.4n}

\begin{definition}\label{2.3n}
 A current $T \in \lnormal_m(M)$ is 
\emph{homologically (mass-)\hspace{0pt}minimiz\-ing} if
 $$
  \mass_U(T) \le \mass_U(T + \rand S)
 $$
 whenever $U \subset M$ is open and relatively compact and 
 $S \in \normal_{m+1}(M)$ has support in $U$.
\end{definition}

For the rest of this subsection we additionally assume that $M$ is compact.
We approach the question of existence of a mass-minimizing current in
every real homology class.
Using de Rham's Theorem, the Hahn--Banach Theorem and the fact that real
Lipschitz $m$-chains are $\rflat$-dense in $\normal_m(M)$, 
cf.\ Section~\ref{s:2.3n}, one can conclude that the homology 
$\Hmmr$ of the chain complex
$\rand \colon \normal_{m+1}(M) \to \normal_m(M)$ is dual to the de Rham
cohomology $\Hmdrm$.
If $\alpha = [T] \in \Hmmr$ is represented by the closed current
$T \in \normal_m(M)$ and $\beta = [\omega] \in \Hmdrm$ is represented
by the closed form $\omega \in \Omega^mM$, then the natural pairing
$\Hmmr \times \Hmdrm \to \R$ is given by
$$
 [\alpha, \beta] = T(\omega) \,.
$$

\begin{proposition}\label{2.4n}
 For every $\alpha \in \Hmmr$ one has
  $$
   \stab(\alpha) = \min \big\{ \mass(T) \bigm|
    T \in \normal_m(M),\ \rand T= 0 \text{ and } [T] = \alpha \big\}\,.
  $$
\end{proposition}

\begin{remark}
In the framework of Geometric Measure Theory it is natural to take the
right hand side of the preceding equation as definition of the stable norm
of~$\alpha$.
We chose the definition of $\stab (\alpha)$ as infimum of the volumes of
Lipschitz cycles representing $\alpha$, since it is geometrically intuitive.
\end{remark}

\begin{proof}
By \cite[3.9]{Fe2} the minimum on the right hand side is attained.
Denote, for the moment, the quantity on the right hand side by 
$\stab'(\alpha)$. 
By Example~1 above, clearly $\stab'(\alpha) \le \stab(\alpha)$ for every
$\alpha \in \Hmmr$.

For the converse, note that since
$\stab$ and $\stab'$ are both norms on $\Hmmr$, it suffices to show
that they coincide on the dense subset $H_m(M,\Q)$.
Let $\alpha \in H_m(M,\Q)$, $T \in \normal_m(M)$ with $\rand T=0$ and
$[T]=\alpha$, and let $\epsilon > 0$. 
We suppose that $M$ is isometrically embedded into some $\R^N$. 
Choose a tubular neighborhood $U$ of $M$ in $\R^N$ so small that
the nearest point projection $p \colon U \to M$ satisfies
$\on{Lip}(p)^m \le \min\big\{1 + \epsilon/\mass(T), 2\big\}$.
By \cite[5.8]{Fe2} there is a closed rectifiable current $R \in \rect_m(M)$
(in fact an integral Lipschitz chain) and $k \in \N$ such that $\frac1k R$ is
homologous to $T$ and $\mass \big(\frac1k R\big) \le \mass (T) + \epsilon$.
(For the notion of rectifiable current see the beginning of
Section~\ref{s:4n}.)

By \cite[Lemma 4.2.19]{Fe} there exists a closed integral polyhedral chain
$P \in \mathcal P_m(U)$ 
(i.e.\ a linear combination of affine simplices with integer coefficients) 
and a $C^1$-diffeomorphism $f$ of $U$ such that 
$f_\# P$ is homologous to $R$ in $U$ and 
$\mass(f_\# P - R) \le \epsilon$. 
After a suitable simplicial subdivision we can write $P$ in the form
$P = \sum_{i=1}^l n_i \current{\Delta_i}$ where the $\Delta_i$ are affine
simplices belonging to a simplicial complex.

Therefore and since $f$ is a diffeomorphism, 
when calculating the mass there is no cancelation.
Putting 
$\tilde c:= \sum \frac{n_i}k \tilde \sigma_i$, 
where $\tilde \sigma_i = f \restr \Delta_i$, we get
$\vol_m(\tilde c) = \frac1k \mass(f_\# P) \le \mass (\frac1k R) + \epsilon \le
\mass(T) + 2 \epsilon$.

Projecting $\tilde c$ to $M$, we get the desired Lipschitz cycle $c$
representing $\alpha$ and
satisfying
$$
 \vol_m(c) \le \on{Lip}(p)^m \vol_m(\tilde c) \le
\big(1+\epsilon/\mass(T)\big)\mass(T) + 4\epsilon  
 \le \mass(T) + 5 \epsilon \,.
$$
\end{proof}

Note that a
 closed current $T \in \normal_m(M)$ with $[T] = \alpha$
 is 
homologically mass-min\-i\-miz\-ing if and only if
$\mass(T)=\stab(\alpha)$.
Such a current is
called a 
\emph{(mass\nobreakdash-)\hspace{0pt}minimizing current in $\alpha$}.
For $m=1$ and $m=n-1$ the structure of homologically mass-minimizing closed
currents is well understood, see \cite{Ba2} for the case $m=1$ and 
\cite{AB'}, \cite{AB} and Section~\ref{s:5n}
for the case $m=n-1$.

The norm on $\Hmdrm$ induced by the comass norm on $\Omega^m M$ is 
equally called \emph{comass norm}, and will be denoted by
$$
 \stab^*(\beta) := \inf \big\{\norm\omega_\infty \bigm|
 \omega \in \Omega^mM,\ d\omega=0 \text{ and } [\omega]=\beta \big\} \,.
$$
It is known that the comass norm on $\Hmdrm$ is dual to the stable norm
$\stab$ on $\Hmmr$, i.e., for all $\alpha \in \Hmmr$ we have
\begin{equation}\label{e:1n}
 \stab(\alpha) = \sup \big\{ [\alpha,\beta] \bigm|
  \beta \in \Hmdrm, \stab^*(\beta) \le 1 \big\}\,,
\end{equation}
see \cite[4.10]{Fe2} or \cite[4.35]{Gr}.
Since we need the arguments from the proof of this duality we will reprove it 
in Section~\ref{s:3.4n}, Theorem~\ref{3.5n}.

\section[Subderivatives and Calibrations]{Subderivatives of the Stable Norm and Calibrations}\label{s:3n}

\subsection{Calibrations}\label{s:3.1n}

A \emph{flat $m$-cochain} is a linear functional on $\normal_m(M)$ that is
continuous with respect to the flat norm $\rflat$.
It is called a \emph{flat cocycle} if it vanishes on the space of boundaries
$$
 \mrboundm = \big\{T \in \normal_m(M) \bigm| 
    \exists S \in \normal_{m+1}(M) \colon \rand S = T \big\}
$$
Note that for a flat cocycle $L$ its flat norm
$$
 \flt(L) = \sup \big\{L(T) \bigm| T \in \nmm, \flt(T) \le 1 \big\}
$$
coincides with its comass norm
$$
 \mass(L) = \sup\big\{L(T) \bigm|  T \in \nmm, \mass(T) \le 1 \big\} \,.
$$
Every flat $m$-cochain can naturally be extended to the $\flt$-closure
$\rflat_m(M)$ of $\nmm$, the space of $m$-dimensional \emph{flat chains}.
These flat chains appear in Lemma~\ref{3.2n}, where they are needed in
the course of the proof, even 
if one restricts the statements to the case of normal currents.
They will be avoided in the rest of the paper.

\begin{definition}\label{3.1n}
 A flat cocycle $L$ of norm $\flt(L)=1$ is called a \emph{calibration}.
 If $L$ is a calibration and $T \in \fmm$ satisfies $L(T) = \mass(T)$,
 then $L$ is said to \emph{calibrate $T$}.
\end{definition}
If $T \in \nmm$ is  calibrated by a
calibration $L$ then $T$ is homologically minimizing.
Indeed, if $S \in \normal_{m+1}(M)$ then
$$
 \mass(T) = L(T) = L(T + \rand S) \le \mass(T + \rand S) \,.
$$

\begin{lemma}\label{3.2n}
Let $L$ be a calibration.
\begin{enumerate}[(a)]
\item
 If $(T_i)_{i \in \N}$ is a sequence in $\fmm$ that $\flt$-converges to
 $T \in \fmm$ and if $L$ calibrates each $T_i$, then $L$ calibrates $T$.
\item
 If $S \in \fmm$ is a piece of $T \in \fmm$, \ie\ if
 $\mass (T) = \mass(S) + \mass(T-S)$, then
 $L$ calibrates $T$ if and only if $L$ calibrates $S$ and $T-S$.
\item
 Assume $L$ calibrates $T \in \fmm$.
 If
 $g\colon M \to \R_{\ge 0}$ is 
 $\mu_T$-integrable 
 then $L$ calibrates $T \inner g$ and 
 $L(T \inner g) = \mass(T \inner g) = \int g\, d\mu_T$.
 If $h \in L^1(M,\mu_T)$ 
 then $L(T \inner h) = \int h\,d\mu_T$.
\end{enumerate}
\end{lemma}
\begin{proof}
Since the mass is lower semicontinuous with respect to flat convergence
we have
$$
 \mass(T) \le \liminf_{i \to \infty} \mass (T_i)
    = \lim_{i \to \infty} L(T_i) = L(T) \,.
$$
This proves (a). 
Statement (b) follows directly from the definitions.
This proves also (c) for the case of step functions.
Using (a) and approximation by step functions one obtains (c).
(Note that the currents $T\inner g_n$, where the $g_n$ are
step functions approximating $g$, need not belong to $\nmm$ even if 
$T\inner g$ does.)
\end{proof}

\begin{definition}\label{3.3n}
 A locally normal $m$-current $T \in \lnmm$ \emph{is calibrated} by the calibration $L$
 if $T \inner A$ is calibrated by $L$ for every compact set $A$ such that
 $T \inner A \in \nmm$.
\end{definition}

\subsection{Subderivatives of the stable norm}\label{s:3.2n}

In this subsection we assume that $M$ is compact.
A subderivative of the stable norm 
$\stab\colon \Hmmr \to \R$ at $\alpha \in \Hmmr$
is a linear form $l \in \Hmmr^*$ such that 
$l(\alpha) = \stab(\alpha)$ and 
$l(\beta) \le \stab(\beta)$ for all $\beta \in \Hmmr$.

\begin{lemma}\label{3.4n}
 Let $l \in \Hmmr^*$ be a subderivative of $\stab$ at $\alpha \in \Hmmr$.
 Then there exists a calibration $L \in \nmm^*$ such that 
 $L(T) = l\big([T]\big)$ for every closed $T \in \nmm$.
 In particular, such an $L$ calibrates every minimizing $T \in \nmm$ in the
 homology class $\alpha$.
\end{lemma}

\begin{remark}
 Since $\stab$ is convex there exists a subderivative $l$ of $\stab$ at 
 $\alpha$ for every $\alpha \in \Hmmr$.
 Hence Lemma~\ref{3.4n} implies that every closed $T \in \nmm$ that
 minimizes mass in $[T]$ is calibrated by some calibration $L$.
\end{remark}

\begin{proof}
For closed currents $T \in \nmm$ we define $L(T) = l\big([T]\big)$.
Since $L(T) = l\big([T]\big) \le \stab\big([T]\big) \le \mass(T)$,
we can use the Hahn--Banach Theorem to extend $L$ to a linear functional
on all of $\nmm$ such that $\mass(L) \le 1$.
Since $L$ vanishes on $\mrboundm$, $L$ is indeed a calibration.
If $T \in \nmm$ is closed,
$[T] = \alpha$ and $\mass(T) = \stab(\alpha)$, then
$$
 L(T) = l(\alpha) = \stab (\alpha) = \mass(T)\,,
$$
i.e., $L$ calibrates $T$.
\end{proof}

\subsection{The canonical representative of a flat cochain}\label{s:3.3n}

According to their definition flat $m$-cochains are objects purely from
functional analysis.
They are elements of $\nmm^*$ that are continuous with respect to the flat
norm.
But it is well known that flat $m$-cochains can be represented by bounded
Lebesgue measurable $m$-forms in the following sense,
cf.\ \cite[IX, Theorem 5A]{Wh} or \cite[4.1.19]{Fe}.

If $L$ is a flat $m$-cochain 
 then, by \cite[IX, Theorem 5A]{Wh} or
\cite[4.1.19]{Fe}
there exists a bounded Lebesgue-measurable $m$-form $\lambda$ such that
for every  smooth $m$-current 
$T_\eta$, 
$\eta \in \Omega^{n-m}_0 M$,
we have
\begin{equation}\label{e:3.3(1)} 
  L(T_\eta) = \int_M \eta \wedge \lambda \,.
\end{equation}
We say that \emph{$\lambda$ is a representative of $L$} 
or that \emph{$\lambda$ represents $L$}.

If $L$ is closed, \ie\ a cocycle, then $\lambda$ is weakly closed, \ie,
$\int_M d\theta \wedge \lambda = 0 $
for every $\theta \in \Omega^{n-m-1}M$,
and we have 
$\rflat (L) = \esssup_{x \in M} \norm {\lambda(x)}$,
where $\norm\ $ denotes the comass norm 
(which coincides with
the Euclidean norm if $m=1$ or $m=n-1$), cf.~Section~\ref{s:2.2n}.

For the proof of Theorem~\ref{reg-reg}
it will be important that
we can find a canonical representative of $L$, 
denoted by $D_L$, by a process of differentiation.
For the following discussion, based on Whitney's book \cite{Wh}, 
we work in standard Euclidean space $\R^n$.
Using charts we can apply the results to manifolds (see below).

Given an oriented $m$-dimensional affine simplex $\sigma$, 
we denote by $P(\sigma)$
the oriented $m$-dimensional affine subspace of $\R^n$ that contains $\sigma$,
by $\xi(\sigma) \in \Lambda_m \R^n$ the unit $m$-vector orienting 
$P(\sigma)$, 
and by $\vol(\sigma) = \mass\big(\current \sigma\big)$ 
the $m$-dimensional Euclidean volume of $\sigma$.

The \emph{thickness} (or \emph{fullness}) $\Theta(\sigma)$ of $\sigma$ is
defined by
$$
  \Theta(\sigma) = \frac{\vol(\sigma)}{\on{diam} (\sigma)^m} \,.
$$
If $p \in \R^n$, 
then the \emph{$p$-thickness} $\Theta_p(\sigma)$ of $\sigma$ is defined by
$$
\Theta_p(\sigma) = \frac{\vol(\sigma)}{\on{diam} (\sigma \cup \{p\})^m} \,.
$$

The following definition is implicit in \cite[IX, \S\,4]{Wh}.

\begin{definition}\label{Wh-regular}
Let $L$ be a flat cochain of degree $m$ on 
$\R^n$ and $p \in \R^n$.
We say that $p$ is a \emph{regular point} for $L$ if there exists an 
$m$-covector $\phi \in \Lambda^m \R^n$ such that 
for every $m$-vector $\xi \in \Lambda_m\R^n$ and all $\epsilon, \eta > 0$
 there exists  $\delta > 0$ with the
property that every $m$-simplex
$\sigma \subset \R^n$  with $p \in \sigma$, $\xi(\sigma) = \xi$,
$\on{diam}(\sigma) < \delta$, and $\Theta(\sigma) \ge \eta$ satisfies
$$
 \left| \langle \xi, \phi \rangle - 
        \frac{L\big(\current \sigma \big)}{\vol(\sigma)} \right |
 < \epsilon \,.
$$
Clearly in this case the covector $\phi$ is unique. We denote it by $D_L(p)$.

The set of all regular points for $L$ is denoted by $\reg L$, its
complement, the set of all \emph{singular points}, by $\sing L$.
\end{definition}
By \cite[IX, Theorem 5A]{Wh}, 
$\sing L$
is a Lebesgue null set and the function 
$D_L \colon \R^n\setminus \sing L \to \Lambda^m \R^n$ is
measurable.
Moreover, $\norm {D_L(p)} \le \rflat(L)$
and the bounded $m$-form $\lambda = D_L$ represents $L$ in the sense 
of~\eqref{e:3.3(1)}

By \cite[X, Theorem 9A]{Wh}, regularity is invariant under 
local diffeomorphisms and $D_L(p)$ behaves like an $m$-covector. 
Hence, using local charts, $\reg L$ and $D_L$ are also defined for 
flat cochains $L$ on a Riemannian manifold.

Proposition~\ref{blowup cochain} below proves that $D_L(p)$ can be 
obtained by a local blow-up of $L$ at $p$.

For $p \in \R^n$ and $r\ge 0$ we denote by
$\mu_{r,p} \colon \R^n \to \R^n$ the homothety with center $p$ and factor $r$,
i.e.,
$\mu_{r,p}(p+x) = p + rx$.

For $L \in \nmrn^*$ a flat $m$-cochain, $p \in \R^n$, and $r>0$, we 
set
$L_{r,p} := \frac1{r^m} (\mu_{r,p})^\# L$, i.e.,
for every $T \in \nmrn$ we have
$L_{r,p}(T) = \frac1{r^m} L\big((\mu_{r,p})_\# T\big)$.
It is easy to see that if $L$ is represented
by the $m$-form $\lambda$ then $L_{r,p}$ is represented by
$\lambda_{r,p}$ where $\lambda_{r,p}(p+x) = \lambda(p +  r x)$.

If $p \in \reg L$, let $L_{0,p}$ denote the flat $m$-cochain 
represented by the constant
$m$-form $\lambda_{0,p}$, where $\lambda_{0,p}(x) := D_L(p)$.

\begin{proposition}\label{blowup cochain}
Let $L$ be a flat  $m$-cochain on
$\R^n$.
\begin{enumerate}[(a)]
\item
 If $p \in \R^n$ and $0 < r \le 1$ then $\rflat(L_{r,p}) \le \rflat(L)$.
\item
 If $p \in \reg L$ then, for $r \to 0$, 
 $L_{r,p}$ converges to $L_{0,p}$ in the weak-$*$-topology on 
 $\normal_m(\R^n)^*$,
 i.e.,
 $\lim\limits_{r\to 0} L_{r,p}(T) = L_{0,p}(T)$ for every 
 $T \in \nmrn$.
\end{enumerate}
\end{proposition}

\begin{proof}
We may assume that $p = 0$.
In the following we omit the subscript $p$.
By \cite[4.1.14]{Fe}, for $0\le r \le 1$ and every  $T \in \nmrn$
we have
$$
 \rflat \big((\mu_r)_\# T \big) \le \max\{r^m,r^{m+1}\} \rflat(T) \le 
 r^m \flt(T).
$$
Therefore,
\begin{align*}
\bigl|L_r(T)\bigr| &= \frac1{r^m} \Big|L\big((\mu_r)_\# T\big) \Big|
 \le \frac1{r^m} \rflat(L) \, \rflat \big((\mu_r)_\# T \big) 
 \le \rflat(L) \, \rflat (T)\,.
\end{align*}
This proves (a).
For the proof of (b),
first assume $T= \current \sigma$ is given by an oriented 
$m$-simplex $\sigma$.
Denoting $\sigma_r := \mu_r \sigma$ we have
$$
 L_r(T) = \frac1{r^m} L\big(\current{\sigma_r}\big) =
 \vol(\sigma) \frac{L\big(\current{\sigma_r}\big)}{\vol(\sigma_r)}
$$
Since
$\Theta_p(\sigma_r) = \Theta_p(\sigma)$ for every $r$,
 \cite[IX, Theorem 10A]{Wh} implies
$$
\lim_{r \to 0}\frac{L\big(\current{\sigma_r}\big)}{\vol(\sigma_r)}
= \big\langle \xi(\sigma), D_L(p) \big\rangle
$$
and thus
$$
\lim_{r \to 0}L_r(T) = \vol(\sigma)\, \big\langle \xi(\sigma), D_L(p) \big\rangle
= L_0(T) \,.
$$
The fact that polyhedral currents are $\rflat$-dense in $\nmrn$ and (a)
imply that
$\lim\limits_{r \to 0} L_r(T) = L_0(T)$ for every $T \in \nmrn$.
\end{proof}

\subsection{Smoothing flat cochains}\label{s:3.4n}

In the proof of Theorem~\ref{theorem} we have to pass from the merely bounded 
and measurable form $\lambda$ representing the flat cochain $L$ to
a smooth approximation, obtained by convolution.
The following lemma is formulated and proved for the case of $M = \R^n$ and
the usual convolution in $\R^n$. 
In the general case, one can embed
$M$ into some~$\R^N$ and perform the convolution using a tubular neighborhood 
of the submanifold~$M$.  See also \cite[4.7]{Fe2}.

\begin{lemma}\label{l:smoothing}
Suppose $L \in \nmrn^*$ is a flat cochain, 
represented by the bounded, measurable $m$-form $\lambda$.
For $\epsilon > 0$, let $\lambda_\epsilon$ originate from $\lambda$ by 
convolution with the kernel $\phi_\epsilon$. 

Then we have
$$
 \lim_{\epsilon \to 0} T(\lambda_\epsilon) = L(T)
$$
for every $T \in \nmrn$.
\end{lemma}

\begin{remark}
 Note that the statement is standard if $\lambda$ is a smooth form or if
 $T$ is given by a smooth form.
\end{remark}

\begin{proof}
 The mollified currents $T_\epsilon$ are defined by 
 $T_\epsilon (\omega) = T(\omega_\epsilon)$ for every 
 $\omega \in \Omega^m_0 \R^n$.
 $T_\epsilon$ is a smooth current, 
 $T_\epsilon = T_{\eta_\epsilon}$ for some 
 $\eta_\epsilon \in \Omega_0^{n-m}\R^n$.
 In particular, $T_\epsilon(\lambda) = \int \eta_\epsilon \wedge \lambda$ 
 is defined
 and $T_\epsilon(\lambda) = L(T_\epsilon)$.
 For $\epsilon \to 0$, $T_\epsilon$ converges to $T$ in the flat norm, and,
 since $L$ is continuous \wrt\ $\flt$-convergence, we have
 $\lim_{\epsilon \to 0} L(T_\epsilon) = L(T)$.
 
 Hence, 
 it suffices to show that $T_\epsilon(\lambda) = T(\lambda_\epsilon)$ holds
 also for the merely bounded and Lebesgue-measurable form $\lambda$.

 Since the support of $T$ is compact, we may assume that also the support
 of $\lambda$ is compact and hence that $\lambda$ is integrable \wrt\
 Lebesgue measure. 
 So we can find a sequence of smooth forms $\lambda^i$ converging to $\lambda$
 in $L^1$.
 Let $\lambda^i{}_\epsilon$ originate from $\lambda^i$ by convolution with
 $\phi_\epsilon$. 
 Then
 $$
  \lambda^i{}_\epsilon(x) - \lambda_\epsilon(x) =
  \int \phi_\epsilon(x-y) \big(\lambda^i(y) - \lambda(y) \big) \, dy
 $$
 and hence
 $$
  \big|\lambda^i{}_\epsilon(x) - \lambda_\epsilon(x)\big|
  \le \sup \betrag{\phi_\epsilon} \,
  \int\big| \lambda^i(y) - \lambda(y) \big| \, dy \,.
 $$
 Therefore, $\lambda^i{}_\epsilon$ converges, for $i \to \infty$, 
 uniformly to $\lambda_\epsilon$.
 So for the current $T = \mu_T \inner \vec T$, we have
 $$
  \lim_{i \to \infty} T(\lambda^i{}_\epsilon) 
  = \lim_{i \to \infty} 
         \int \langle \vec T, \lambda^i{}_\epsilon \rangle \, d\mu_T
  = \int \langle \vec T, \lambda_\epsilon \rangle \, d\mu_T
  = T(\lambda_\epsilon) \,.
 $$
 On the other hand, since $\lambda^i$ is smooth, 
 we have $T(\lambda^i{}_\epsilon) = T_\epsilon(\lambda^i)$, 
 and since $T_\epsilon$ is smooth, $L^1$-convergence of 
 $\lambda^i \to \lambda$ yields
 $$
  \lim _{i \to \infty} T(\lambda^i{}_\epsilon) 
  = \lim _{i \to \infty} T_\epsilon(\lambda^i) = T_\epsilon(\lambda) \,.
 $$
 So $T(\lambda_\epsilon) = T_\epsilon(\lambda) = L(T_\epsilon)$
 and 
 $\lim\limits_{\epsilon \to 0} T(\lambda_\epsilon)
  = \lim\limits_{\epsilon \to 0} L(T_\epsilon) = L(T)$. 
\end{proof}

Using Lemma~\ref{3.4n} and Lemma~\ref{l:smoothing} 
we can easily prove~\eqref{e:1n}:

\begin{theorem}\label{3.5n}
 On a compact and oriented Riemannian manifold $M$ the comass norm
 $\stab^*$ on $\Hmdrm$ is dual to the stable norm $\stab$ on $\Hmmr$
 with respect to the natural pairing between homology and cohomology.
\end{theorem}

\begin{remark}
 This statement follows from \cite[4.10]{Fe2}. The proof given here
 elaborates the one sketched in \cite[4.35]{Gr}.
\end{remark}

\begin{proof}
 If $\alpha \in \Hmmr$ is represented by $T \in \nmm$ and if 
 $\beta \in \Hmdrm$ is represented by $\omega \in \Omega^mM$ then
$$
 [\alpha, \beta] = T(\omega) \le \mass(T) \norm \omega_\infty \,.
$$
This implies 
$$
 [\alpha, \beta] \le \stab(\alpha) \, \stab^*(\beta) \,.
$$
It remains to show that
\begin{equation}\label{e:2n}
 \stab(\alpha) \le \sup_{\stab^* (\beta) \le 1} [\alpha, \beta] \,.
\end{equation}
 According to Proposition~\ref{2.4n} we can choose a minimizing 
 $T \in \nmm$ in $\alpha$
 and according to Lemma~\ref{3.4n} there exists a calibration $L$ 
 calibrating $T$.
 Then
 $$
  L(T) = \mass(T) = \stab(\alpha) \,.
 $$
 By convolution we mollify a measurable $m$-form $\lambda$ representing $L$
 to obtain closed forms $\lambda_\epsilon \in \Omega^mM$ such that
 $$
  \lim_{\epsilon \to 0} \norm{\lambda_\epsilon}_\infty =
  \esssup_{x \in M} \norm{\lambda(x)} = 1 \,.
 $$
 Now Lemma~\ref{l:smoothing} implies 
 $$
  \stab(\alpha) = L(T) = \lim_{\epsilon \to 0} T(\lambda_\epsilon)
  = \lim_{\epsilon \to 0} [\alpha, \beta_\epsilon]
 $$
 where $\beta_\epsilon = [\lambda_\epsilon] \in \Hmdrm$ satisfies
 $$
  \liminf_{\epsilon \to 0} \stab^*(\beta_\epsilon) \le 
  \lim_{\epsilon \to 0} \norm{\lambda_\epsilon}_\infty = 1 \,.
 $$
 This proves \eqref{e:2n}.
\end{proof}

\section[Calibrations in Codimension One]
{Calibrations and Minimizing Currents in Codimension One}
\label{s:4n}

In this section we will show that the singular set of a codimension one
calibration $L$ contains the singular sets of all closed integer multiplicity
rectifiable currents calibrated by $L$. 
For the definition of the 
set $\hlrect(M)$ 
of integer multiplicity rectifiable
currents, for the definition of the regular part $\reg T$ of a 
homologically minimizing $T \in \hlrect(M)$ and for the regularity theory
for such currents we refer to \cite[\S 27 and \S 37]{Si}.
Here, we 
note that if $T \in \hlrect(M)$  
and if $\rand T$ has locally finite mass, 
then $T \in \hlnormal(M)$.
Moreover, if $T \in \hlrect(M)$
is closed and  
homologically minimizing then
$\reg T$ is a smooth hypersurface in $M$ oriented by a smooth unit
$(n{-}1)$-vector field $\vec T$ and 
$$
   T(\omega) = \int_{\reg T} \omega
$$
for all $\omega \in \Omega_0^{n-1}M$.

\begin{lemma}\label{reg-det}
Suppose that the flat cocycle $L \in \hnormal(M)^*$ is a calibration and that 
 $T \in \hlrect(M)$ is
closed and
calibrated by $L$.

 Then, for every $p \in \reg T \cap \reg L$, 
 we have $D_L(p) = \vec T(p)^\flat$.
\end{lemma}

Here, $\vec T(p)^\flat \in \Lambda^{n-1} T_p M$ denotes the $(n{-}1)$-covector satisfying
$\langle  \xi , \vec T(p)^\flat \rangle = g_p ( \xi, \vec T(p))$
for all $(n{-}1)$-vectors $\xi \in \Lambda_{n-1} T_p M$. 
For the definition of $\vec T$ see Section~\ref{2.2n}.

\begin{proof}
Let $p \in \reg T \cap \reg L$. 
Since the statement is local we may work in a local chart, \ie\
in 
$\R^n$ equipped with a Riemannian metric $g$.
Mass and volume will be defined \wrt\ $g$.
To distinguish metric terms which refer to $g$ from the Euclidean ones, we
will mark them by a subscript or superscript $g$.
We may assume that $p$ is the origin, that the metric $g_p$ at the origin
coincides with the standard Euclidean scalar product and that
$\reg T$ is a hyperplane $P$ through the origin, oriented by the unit 
$(n{-}1)$-vector $\xi:= \vec T(p)$.

For any $(n{-}1)$-simplex $\sigma$ in $P$ with the orientation of $T$, 
we have
$L\big(\current \sigma \big) = \mass_g\big(\current \sigma \big)
= \vol_g(\sigma)$.
Therefore, every sequence $\sigma_i$ of $(n{-}1)$-simplices with
$p \in \sigma_i$ and $\xi(\sigma_i) = \xi$ for every $i \in \N$,
and $\lim_{i \to \infty} \diam (\sigma_i) = 0$ satisfies
$$
  \lim_{i \to \infty} \frac{L\big(\current \sigma_i \big)}{\vol (\sigma_i)}
 =\lim_{i \to \infty} \frac{\vol_g(\sigma_i)}{\vol (\sigma_i)}
  = 1 \,.
$$
Recalling Definition~\ref{Wh-regular} we get $\langle \xi,D_L(p) \rangle = 1$. 
Since $\betrag{D_L(p)}= \norm {D_L(p)} \le 1$ and since $\betrag \xi = 1$, 
this implies that
$D_L(p) = \xi^\flat = \vec T(p)^\flat$.
The crucial point is that the comass norm on the space of 
$(n{-}1)$-covectors is Euclidean and hence strictly convex.
\end{proof}

\begin{theorem}\label{reg-reg}
 Suppose that the flat cocycle $L \in \hnormal(M)^*$ is a calibration and that 
 $T \in \hlrect(M)$ is closed and calibrated by $L$.
 Then $\sing T \subset \sing L$. 

 In particular, the union of the singular sets of all closed 
 $T \in \hlrect(M)$ 
 calibrated by $L$ is a Lebesgue null set.
\end{theorem}

\begin{proof}
Note that $T$ is homologically mass-minimizing. 
Assume that $p \in \spt T$ is a regular point of $L$. 
We have to show that $p \in \reg T$.
As above, since the statement is of local nature, we may work in local
coordinates. 
By the Decomposition Theorem for codimension one rectifiable currents,
cf.~\cite[27.6]{Si},
we may assume that $T$ is of multiplicity one.
We identify some neighborhood of $p$ in $M$ with 
$\R^n$, equipped with a Riemannian metric $g$, the point $p$ with
the origin $0$, 
such that $g_p = g_0$ corresponds to the standard scalar product of $\R^n$.
Mass w.r.t.\ $g$ is denoted by $\mass_g$ while all other metric terms refer
to the Euclidean metric.

The regularity theory for mass-minimizing currents implies that there
exists a tangent cone $C$ of $T$ at $p$, \ie,
there exists a sequence $r_i > 0$ converging to zero such that
the sequence $T_i:=\push{\mu_{1/r_i}} T$ 
converges weakly to a closed multiplicity one 
current $C \in \hlrect(\R^n)$ (cf. \cite[Theorem 37.4]{Si}).
This current $C$ is a cone, \ie, $\push{\mu_r} C = C$ for every $r > 0$. 
The tangent cone $C$ is mass-minimizing with respect to the metric $g_0$.
Hence $C$ is given by the smooth hypersurface $\reg C$ (which is a cone, too)
with a possible singular set of dimension $\le n-8$.

We will show that $C$ is calibrated by the cochain $L_0$ represented by
the constant $(n{-}1)$-form $\lambda_0(x) := D_L(p)$.
From this we will conclude that $\spt C$ is a hyperplane, 
and then regularity theory  shows that $p$ is a regular point of~$T$.

Let $W \subset \R^n$ be an open set with compact closure that
contains the origin such that 
$\tilde C := C \inner W \in \hnormal(\R^n)$.
By the arguments used in the proofs of \cite[Theorem 37.2]{Si} and
\cite[5.4.2]{Fe}, there exists a sequence of compact sets $K_i \subset \R^n$ 
such that
the currents $S_i := T_i \inner K_i \in \hnormal(\R^n)$ satisfy
\begin{align*}
 \lim_{i \to \infty} \rflat(S_i - \tilde C) = 0 \, \quad
\text{and} \quad
 \lim_{i \to \infty} \mass_{g_i}(S_i) = \mass(\tilde C) \,,
\end{align*}
where $g_i := \frac1{r_i^2}(\mu_{r_i})^* g$. 
Note that $S_i = (\mu_{1/r_i})_\# (T \inner \tilde K_i)$ where
$\tilde K_i := \mu_{r_i} K_i$.

If we set $L_i := \frac{1}{r_i^{n-1}}( \mu_{r_i})^\# L$, 
 Proposition \ref{blowup cochain} (b)
gives us $\lim\limits_{i \to \infty}L_i(S) = L_0(S)$ 
for every $S \in \hnormal(\R^n)$.
By Proposition \ref{blowup cochain} (a) the sequence of cochains
$L_i$ is uniformly bounded w.r.t.\ the flat norm. 
Hence
$$
 \lim_{i \to \infty}L_i(S_i) = L_0 (\tilde C) \,,
$$
where $L_0$ is given by the constant $(n{-}1)$-form $\lambda_0(x) = D_L(p)$.
Since $L_0$ is closed and $\rflat(L_0) = \norm {D_L(p)} \le 1$,
$L_0$ is a calibration.

Using the fact that $L$ calibrates $T$, we get
\begin{align*}
 L_i(S_i) & = \frac{1}{r_i^{n-1}} L\big((\mu_{r_i})_\# S_i \big)
 =  \frac{1}{r_i^{n-1}} L\big((\mu_{r_i})_\# (\mu_{1/r_i})_\# (T \inner \tilde K_i) \big) \\
& = \frac{1}{r_i^{n-1}} L (T \inner \tilde K_i)  = 
\frac{1}{r_i^{n-1}} \mass (T \inner \tilde K_i) = 
\mass_{g_i}( S_i) \,.
\end{align*}
Therefore, 
$L_0(\tilde C) = \lim\limits_{i \to \infty}  L_i(S_i) = 
 \lim\limits_{i \to \infty} \mass_{g_i}( S_i) = \mass(\tilde C)$,
i.e., $L_0$ calibrates~$\tilde C$.

Now, Lemma~\ref{reg-det} implies that $\vec C(x)^\flat = D_L(p)$ for
every $x \in \reg C \cap W$.
In particular, $\vec C(x)$ does not depend on 
$x \in \reg C$. 
Since $\rand C=0$, we can conclude that $\spt C$ is a hyperplane.
Since $C$ is of multiplicity one this implies that the density of $C$ at $p$,
and hence the density of $T$ at $p$, is one, 
and it follows from
\cite[5.4.6]{Fe} together with \cite[5.4.5 (2)]{Fe}
that $p \in \reg T$.
This proves the first statement.

The second one follows from the fact that, by \cite[IX, Theorem 5A]{Wh},
$\sing L$ is a Lebesgue null set.
\end{proof}

\section[Minimizing Closed Currents in
          Codimension One]{The Structure of Homologically Minimizing Closed Currents in
          Codimension One} \label{s:5n}

\subsection{The case of boundaries} \label{s:5.1n}

It is not difficult to prove that $\lnormal_n(M)$ is precisely the set
of currents of the form $T_f = \current M \inner f$, where 
$f \in \lbv(M)$ is of locally bounded variation,
cf.\ \cite[4.5.7]{Fe},
and 
$$
  T_f(\omega) = \int_M f \, \omega \,
$$
for all $\omega \in \Omega_0^{n}M$.
Hence a codimension one boundary $T \in \hlnormal(M)$ is given as
$T = \rand T_f $ for some
$f \in \lbv(M)$.
Additionally, we will assume that $T$ is homologically minimizing,
cf.\ Definition~\ref{2.3n}.
We will give an overview over results on the structure of such $T$,
which are proved in \cite[4.5.9]{Fe} and \cite{AB'}, \cite{AB}.

The BV-function $f$ can be chosen to be (upper or lower) semicontinuous.
For $s \in \R$ consider the sets
$\{x \in M \mid f(x) > s \}$ and $\{x \in M \mid f(x) \ge s \}$ and
let
$T_{s+} := \rand \current{\{x \in M \mid f(x) > s \}}$ and
$T_{s-} := \rand \current{\{x \in M \mid f(x) \ge s \}}$.
Here, if $A \subset M$ is measurable then $\current A \in (\Omega_0^n M)^*$
denotes the current defined by $\current A (\omega) = \int_A \omega$
for $\omega \in \Omega_0^n M$.
We have:
\begin{enumerate}
\item
 $T_{s+} = T_{s-}$ for all but countably many $s \in \R$.
\item
 $T_{s+}, T_{s-} \in \hlrect(M)$ for every $s \in \R$.
\item
 Every $T_{s+}, T_{s-}$ is homologically minimizing, so
\item \label{l:4}
 for each $T_{s\pm}$, its regular part $\reg T_{s\pm}$ 
 is a smooth hypersurface,
 dense in $\spt T_{s\pm}$, and the singular part 
 $\sing T_{s\pm}=\spt T_{s\pm} \setminus \reg T_{s\pm}$ is of
 Hausdorff dimension at most $n-8$.
\item
 $T = \int_\R T_{s+}\,ds = \int_\R T_{s-}\,ds$.
\item\label{l:6}
 $\mass_g(T) = \int_\R \mass_g(T_{s+})\,ds = \int_\R \mass_g(T_{s-})\,ds$.
\item\label{l:7}
 $\spt T = \bigcup_{s \in \R}\big( \spt T_{s+} \cup \spt T_{s-}\big)$.
\item\label{l:8}
 For every $s \in \R$,
 $$
  T_{s+} = \operatorname*{\rflat-lim}_{h \to 0+} \frac1h
            \int_s^{s+h} T_{t\pm}\, dt
  \quad \text{and} \quad 
  T_{s-} = \operatorname*{\rflat-lim}_{h \to 0+} \frac1h
            \int_{s-h}^{s} T_{t\pm}\, dt \,.
 $$
\end{enumerate}
 Let $J = \{s \in \R \mid T_{s-} \ne T_{s+} \}$.
\begin{definition}\label{def-gap}
For $s \in J$ let $G_s$ denote the interior of  $\{x \in M \mid f(x) = s \}$.
The sets $G_s$ are called the \emph{gaps} of $T$.
\end{definition}
Then we have
\begin{align*}
 M \setminus \spt T 
  &= \bigcup _{s\in J} G_s \,.
\end{align*}

For the boundary of $G_s$ in the sense of currents, we have
$\rand \current{G_s} = T_{s-} - T_{s+}$.
The set theoretic boundary $\rand G$ of $G_s$ 
consists of connected components of $\spt T_{s-} \cup \spt T_{s+}$.
Outside the singular sets $\sing T_{s-} \cup \sing T_{s+}$ 
this boundary is smooth and
consists of  connected components of 
the hypersurfaces $\reg  T_{s-}$ and  $\reg T_{s+}$.
We denote the smooth part of $\rand G$ by $\reg \rand G$.

We call  $\reg T := \bigcup_{s \in \R}\big(\reg T_{s+} \cup \reg T_{s-}\big)$
the \emph{regular} set of $T$
 and
$\sing T := \bigcup_{s \in \R}\big(\sing T_{s+} \cup \sing T_{s-}\big) =
 \spt T \setminus \reg T$
the \emph{singular} set of $T$.
Although
the regularity theory for mass-minimizing rectifiable $(n{-}1)$-currents
implies that
 the singular set of each $T_{s\pm}$ is of Hausdorff dimension 
$\le n-8$, \emph{a priori} it is not clear whether their union $\sing T$ is
small. However, the following lemma implies
that it is at least a Lebesgue null set, cf.\ Corollary~\ref{sing-null}.

\begin{lemma}\label{determined}
 Suppose $T \in \hlnormal(M)$ is a boundary calibrated by the
 flat cocycle $L$.
 Then:
 \begin{enumerate}[(a)]
 \item
 For every $p \in \reg T \cap \reg L$, the $(n{-}1)$-covector $D_L(p)$ 
  is uniquely determined by~$T$,
  $D_L(p) = \vec T_{s\pm}(p)^\flat$ if $p \in \reg T_{s\pm}$.
 \item
 $\sing T \subset \sing L$. 
 In particular,
 $\sing T$ is a Lebesgue null set.
 \end{enumerate}
\end{lemma}

\begin{proof}
Since $T$ is calibrated by $L$, it is homologically mass-minimizing. 
So we can apply 
the list of statements above.
Points \ref{l:6} and~\ref{l:8} imply that $L$ calibrates each 
$T_{s+}$ and $T_{s-}$, 
cf.\ Lemma~\ref{3.2n}.

Then, by  Lemma~\ref{reg-det}, we have $D_L(p) = \vec T_{s\pm}(p)^\flat$
for every $p \in \reg T_{s\pm}\cap \reg L$.
This proves (a).
Statement (b) follows immediately from Theorem~\ref{reg-reg}. 
\end{proof}

{\sloppy
\begin{corollary}\label{sing-null}
  Suppose $T \in \hlnormal(M)$ 
  is a locally mass-minimizing boundary.
  Then $\sing T$ (as defined above) is a Lebesgue null set.
\end{corollary}
}

\begin{proof}
  Since the statement is of local nature we may assume that $T$ has
  compact support.
  Using the Hahn--Banach Theorem like in Section~\ref{s:3.2n}
  we get a calibration that calibrates $T$ (cf. \cite[4.10]{Fe2}).
  Then the statement follows from Lemma~\ref{determined}.
\end{proof}

\subsection{The covering space $M_\alpha$ associated to 
        $\alpha \in \Hhmr$}\label{s:5.2n}

From now on we assume that $M$ is compact.
In order to apply the results of the preceding subsection to closed
currents that do not bound, we will lift these to an appropriate
(infinite) covering space of $M$.
Here we present the relevant material from topology.

Using currents, we can describe the \emph{Poincar\'e duality} isomorphism
between $\Hidrm$ and $\Hhmr$ as follows (cf.\ \cite[Theorem 14]{dR}).
Every $\alpha \in \Hhmr$ can be represented by a smooth closed current
$T_\eta$, where $\eta \in \Omega^1M$ is closed.
Then the cohomology class of $\eta$ depends only on $\alpha$ and is
called the \emph{Poincar\'e dual} $\alpha^\PD \in \Hidrm$ of $\alpha$.
The natural pairing between $\Himr$ and $\Hidrm$ yields the
\emph{intersection form}
$$
  I \colon \Himr \times \Hhmr \to \R,\ 
  I(h,\alpha) = [h, \alpha^\PD] \,.
$$
Explicitly, if $\alpha = [T_\eta]$ and $h$ is represented by a real Lipschitz
1-cycle $c$, then
\begin{equation}\label{e:3n}
 I(h,\alpha) = \int_c \eta \,. 
\end{equation}

The image of the Hurewicz homomorphism
$H \colon \pi_1(M) \to \Himr$
is the set of all integer classes in $\Himr$ and will be denoted by
$\Himzr$.
With respect to the intersection form $I$ the lattice
$\Himzr$ in $\Himr$ is dual to the lattice $\Hhmz$ in $\Hhmr$.

We set
\begin{equation}\label{e:4n}
K(\alpha) = \big\{ k \in \Himzr \bigm| I(k,\alpha) = 0 \big\}
\end{equation}
and 
$\tilde K(\alpha) = H^{-1}\big(K(\alpha)\big) \subset \pi_1(M)$.
Now the covering
$
 p \colon M_\alpha \to M
$
associated to $\alpha \in \Hhmr$ is given by
$$
  M_\alpha = \tilde M / \tilde K(\alpha)
           = \bar M / K(\alpha) \,,
$$
where $\tilde M$ denotes the universal covering and 
$\bar M = \tilde M / \ker H$
denotes the Abelian covering of $M$.
Equations \eqref{e:3n} and~\eqref{e:4n} imply that 
$p \colon M_\alpha \to M$ is the smallest covering space of $M$ such that
$p^*\eta$ is exact.
The group of deck transformations of $p$ is isomorphic to
\begin{equation}\label{e:5n}
 \Himzr / K(\alpha) \simeq \pi_1(M)/\tilde K(\alpha) \simeq
 \Z^{b-\on{rk} K(\alpha)}
\end{equation}
where $b=b_1(M) = \dim \Himr$ is the first Betti number of $M$.
We denote the deck transformation corresponding to $k \in \Himzr$
by $\tau_k \colon M_\alpha \to M_\alpha$.
Note that $\tau_k \ne \id_{M_\alpha}$ iff $I(k,\alpha) \ne 0$.

Recall that in the introduction we defined $V(\alpha)$ to be the smallest
linear subspace of $\Hhmr$ that is spanned by integer classes and 
contains $\alpha$.
From the preceding discussion 
one can conclude
that $V(\alpha)$ is the orthogonal complement
of $K(\alpha)$ with respect to $I$,
$$
  V(\alpha) = \big\{ \beta \in \Hhmr \bigm| 
         I(k,\beta) = 0 \text{ for all } k \in K(\alpha) \big\} \,.
$$
If $\beta \in V(\alpha)$ then obviously $V(\beta) \subset V(\alpha)$.
Hence the preceding equation implies that $K(\alpha) \subset K(\beta)$
and $\tilde K(\alpha) \subset \tilde K(\beta)$.
This proves:
\begin{itemize} \refstepcounter{equation}\label{e:6n}
\item[(\theequation)\ \ ]
 If $\beta \in V(\alpha)$ and if $\eta \in \Omega^1M$ represents
       $\beta^\PD \in \Hidrm$, then $p^*\eta \in \Omega^1 M_\alpha$ is
       exact. 
\end{itemize}
Now let $\eta \in \Omega^1M$ represent $\alpha^\PD \in \Hidrm$.
Then there exists a primitive $g \in C^\infty(M_\alpha, \R)$ of
$p^* \eta$, i.e., $dg = p^*\eta$, and
\eqref{e:3n} implies
\begin{equation}\label{e:10n}
 g(\tau_k x) = g(x) + I(k,\alpha)
\end{equation}
for all $k \in \Himzr$ and all $x \in M_\alpha$.

\subsection{The lift to $M_\alpha$}\label{s:5.3n}

Let $\alpha \in \Hhmr$ and consider the covering $p \colon M_\alpha \to M$ 
with the induced metric on $M_\alpha$.
The \emph{lift} $p^\#T \in \hlnormal(M_\alpha)$ of $T\in \hnormal(M) $
 to $M_\alpha$ is
defined by
$$
 \big(p^\# T \big)(\omega) = \big(((p \restr U)^{-1})_\# T\big)(\omega)
$$
provided that $\omega \in \Omega^{n-1}M_\alpha$ has compact support in an
open subset $U \subset M_\alpha$ on which $p$ is injective.
Note that $\spt (p^\# T) = p^{-1} (\spt T)$.
For $\eta \in \Omega^1M$ we have $p^\#(T_\eta) = T_{p^* \eta}$.
If $p^*\eta$ is exact, $p^*\eta = dg$, then 
$p^\#(T_\eta)$ is a boundary
$$
   p^\#(T_\eta) = - \rand T_g \,,
$$
cf.\ Section~\ref{s:2.2n}, Example 2.

\begin{lemma}\label{primitive}
Suppose $T \in \hnormal(M)$ is a closed normal current representing 
$\alpha \in \Hhmr$.
Then there exists $f \in \lbv(M_\alpha,\R)$ such that
$$
 p^\# T = \rand T_f
$$
and
\begin{equation}\label{e:quasiper}
  f \circ \tau_k = f - I( k, \alpha )
 \end{equation}
for  every $k \in \Himzr$.
\end{lemma}

\begin{proof}
According to Section~\ref{s:5.2n}, $\alpha$ can also be represented by
a current of the form
$T_\eta$, $\eta \in \Omega^1 M$ a closed smooth 1-form.
Therefore $T = T_ \eta + \rand S$ for some 
$S \in \normal_n(M)$.
$S$ is of the form $S = T_h$ for some $h \in BV(M)$,
cf. Section~\ref{s:5.1n}.
So, if $p^*\eta = d g$ for $g \in C^\infty(M_\alpha,\R)$, then
the lift $p^\# T$ of $T$ to $M_\alpha$ has the form
\begin{align*}
  p^\# T 
   = p^\# (T_ \eta ) +  \rand (p^\# S)
        = - \rand T_{g} 
           + \rand T_{h \circ p}
        = \rand  T_f \,,
\end{align*}
where
$f := - g + h \circ p \in \lbv(M_\alpha)$.
Since $g$ satisfies \eqref{e:10n}, $f$ satisfies~\eqref{e:quasiper}.
\end{proof}

\begin{lemma}\label{lift-cal}
Suppose 
that $T \in \hnormal(M)$ is a closed normal current, that
$L \in \hnormal(M)^*$ is a calibration,
represented by the bounded measurable $(n{-}1)$-form $\lambda$,
and that $L$ 
calibrates $T$.

Then the lift $p^\# L$ of $L$, defined by 
$(p^\# L)(S) = L(p_\# S)$ for every $S \in \hnormal(M_\alpha)$,
is a calibration represented by $p^* \lambda$ 
and calibrates the lift $\bar T:= p^\# T$ of $T$
and every leaf $\bar T_s$ of~$\bar  T$.
\end{lemma}

\begin{proof}
 Since $p$ is distance-nonincreasing we have $\flt(p_\# S) \le \flt(S)$
for every $S \in \hnormal(M_\alpha)$.
This shows that $\flt(p^\# L) \le 1$.
Then the statements follow from Section~\ref{s:3.3n} and Lemma~\ref{3.2n},
and from Points \ref{l:6} and~\ref{l:8} in Section~\ref{s:5.1n}.
\end{proof}

In particular, 
if $T$ is a minimizer in $\alpha$, then its lift $\bar T= p^\# T$ 
to $M_\alpha$ is homologically minimizing.
Therefore it has the properties described
in Section~\ref{s:5.1n} above.
In particular,
$M_\alpha \setminus \spt \bar T$ is the countable union of gaps.
\begin{corollary}\label{gaps}
 For each gap $G$ of $\bar T$ the restriction $p \restr G$ of $p$ to $G$
 is injective. In particular, $G$ has finite volume.
\end{corollary}

\begin{proof}
 By Definition~\ref{def-gap}, $G = G_s$ for some $s \in J \subset \R$.
 If $\tau_k \ne \on{id}_{M_\alpha}$ is a deck transformation of
 $p \colon M_\alpha \to M$ then $I( k, \alpha ) \ne 0$, and 
 \eqref{e:quasiper} implies that $\tau_k G_s \cap G_s = \emptyset$.
 Hence, $p \restr G_s$ is injective
 and $\vol_n(G) \le \vol_n(M) < \infty$.
\end{proof}

\section{Proof of Theorem~\ref{theorem}}\label{s:6n}

It suffices to consider $\alpha \in \Hhmr \setminus\{0\}$.
To show that the restriction $\stab \restr V(\alpha)$ of $\stab$
to $V(\alpha)$ is differentiable at $\alpha$,
we have to prove the following.
If  $l_1$, $l_2$ are two subderivatives of $\stab$ at $\alpha$ and if
$\beta \in V(\alpha)$, then
$l_1(\beta) = l_2(\beta)$.
Let $T \in \hnormal(M)$ be a minimizer in~$\alpha$, \ie, $[T] = \alpha$ and 
$\stab(\alpha) = \mass(T)$.
By Lemma \ref{3.4n} 
there are flat cochains $L_1$ and $L_2$ corresponding to 
$l_1$ and $l_2$, respectively, that calibrate $T$, \ie,
$L_i(S) \le \mass(S)$ for every $S \in \hnormal(M)$ and $L_i(T) = \mass(T)$,
for $i = 1,2$.
By Section~\ref{s:3.3n} there are Lebesgue-measurable 
$(n{-}1)$-forms $\lambda_1:=D_{L_1}$, $\lambda_2:=D_{L_2}$ 
representing the cochains $L_1$, $L_2$ with
$\betrag{\lambda_1} \le 1$, $\betrag{\lambda_2} \le 1$.

We represent $\beta$ by a smooth current $T_\eta$, 
$\eta \in \Omega^1 M$ a smooth
closed 1-form
with $[\eta] = \beta^\PD \in \Hidrm$.
Then 
$l_i(\beta) = L_i(T_\eta) = \int_M \eta \wedge \lambda_i$ for $i = 1,2$.
Using Lemma~\ref{lift-cal} and
applying Lemma~\ref{determined} to $p^\# L$ and $p^\# T$ we conclude that
$\lambda_1 = \lambda_2$ Lebesgue almost everywhere
on $\spt T$. 
Hence it is enough to show that
$$
  \int_{M \setminus \spt T} \eta \wedge (\lambda_1 - \lambda_2) = 0 \,.
$$
Consider the covering $p \colon M_\alpha \to M$ associated to~$\alpha$ and
let $\bar T := p^\# T$ be the lift of $T$ to $M_\alpha$.
Then $M \setminus \spt T = p(M_\alpha \setminus \spt \bar T)$,
where $M_\alpha \setminus \spt \bar T = \bigcup_{s \in J} G_s$ is the union
of at most countably many gaps, cf.\ Section~\ref{s:5n}.
By Corollary~\ref{gaps},
for every gap $G$, the restriction $p \restr G \colon G \to M$ 
of $p$ to~$G$ is injective.
Hence
it suffices to show that
$$
 \int_{G} p^* \big( \eta \wedge (\lambda_1 - \lambda_2) \big) = 0
$$
for every gap $G$.

Since $\beta \in V(\alpha)$, 
the 1-form $p^* \eta$ is exact, cf.~\eqref{e:6n}. 
So there exists
 $g \in C^\infty(M_\alpha,\R)$  such that $p^* \eta = dg$ and 
 we have to show that
$$
\int_{G} dg \wedge (\bar\lambda_1 - \bar\lambda_2)  = 0 \,,
$$ 
where $\bar\lambda_i = p^* \lambda_i$.

Our aim is to cut off the ``ends'' of $G$ and apply Stokes's Theorem to the
remaining compact domain.
For any given $\delta > 0$, we
will show that
$$
\left\vert\int_G dg \wedge \omega\right\vert < 
 \delta\,,
$$ 
where $\omega:= \bar\lambda_1 - \bar\lambda_2$.
Choose $x_0 \in G$ and let $\tilde d$ be a smooth 
function such that 
$\sup_{x\in M_\alpha}\bigl|\tilde d(x) - d(x_0,x)\big| \le 1$ 
and $\on{Lip}(\tilde d) \le 2$.
For $r > 0$ set $B_r := \{x \in M_\alpha \mid \tilde d(x) < r \}$.
Since $p^* \eta = dg$ is bounded,
there exists  $r_0>0$
such that $\sup_{ x \in B_r}\betrag{g(x)} \le 2C r$ for all $r \ge r_0$,
where $C = \sup_{M_\alpha}\betrag{dg}$.

By Corollary~\ref{gaps},
$\vol_n (G)$ is finite. 
By the coarea formula, we have
$$
  \int_0^\infty \vol_{n-1} (G \cap \rand B_r) \, dr \le 
  \on{Lip}(\tilde d) \vol_n(G) < \infty \,.
$$
Therefore we can find $r > r_0$ such that
$\rand B_r$ is a smooth hypersurface which meets the regular part
$\reg \rand G$ of 
$\rand G$ transversely and such that 
$$
\vol_n(G \setminus B_r) < \frac \delta {8C} \text{ and }
\vol_{n-1}(G \cap \rand B_r) < \frac \delta {24Cr}\,.
$$
In particular, since 
$\betrag{\bar\lambda_i} \le 1$, 
we have
$$
 \left|\int_{G \setminus B_r} dg \wedge \omega\right| 
 < \frac\delta4 \,.
$$
Since the $(n{-}1)$-forms $\bar \lambda_1$, $\bar \lambda_2$ are 
not smooth, 
we have to pass to smooth forms in order to apply Stokes's Theorem. 
We do this by mollifying the $\lambda_i$ by means of convolution.
Thus we find smooth approximations $\lambda_i^\epsilon$ of the $\bar \lambda_i$
which satisfy $\sup\betrag{\lambda_i^\epsilon} \le \frac32$ and
$$
   \left|\int_{G \cap B_r} dg \wedge (\omega^\epsilon - \omega) \right|
   < \frac \delta 4 \,,
$$
where $\omega^\epsilon = \lambda_1^\epsilon - \lambda_2^\epsilon$.
Since, for $i = 1,2$, $L_i$ is a cocycle, 
 $\bar \lambda_i$ is weakly closed and hence $d \lambda_i^\epsilon = 0$.
Now, the boundary currents 
$T^+ := T_{s+} \inner \rand G_s$ and $T^- := T_{s-}\inner \rand G_s$ of 
$G = G_s$
are calibrated by the lift $\bar L_i = p^\# L_i$, 
represented by $\bar\lambda_i$,
cf.\ Lemma~\ref{lift-cal}. 
By Lemma~\ref{3.2n} (c)
we have
$$
  \bar L_i \big((T^\pm \inner g) \inner B_r \big) = 
  \int_{\reg T^\pm \cap B_r} g \ d\vol_{n-1} \,.
$$
In particular, 
$\bar L_1 \big((T^\pm \inner g) \inner B_r \big) = 
\bar L_2 \big((T^\pm \inner g) \inner B_r \big)$.
Hence  
Lemma~\ref{l:smoothing} implies
\begin{gather*}
  \lim_{\epsilon \to 0} 
  \left( \int_{\reg T^\pm \cap B_r}  g \lambda_1^\epsilon
       - \int_{\reg T^\pm \cap B_r}  g \lambda_2^\epsilon  \right) = \\ 
    \bar L_1 \big((T^\pm \inner g) \inner B_r \big) 
  - \bar L_2 \big((T^\pm \inner g) \inner B_r \big)
  = 0 \,.
\end{gather*}
So we can choose $\epsilon$ such that
\begin{equation*}
  \left| \int_{\rand G \cap B_r} g \,
          \omega^\epsilon \right| 
  < \frac \delta 4 \,.
\end{equation*}
The inequalities 
$\sup_{B_r} \betrag g \le 2Cr$, $\sup \betrag{\omega^\epsilon} \le 3$
and $\vol_{n-1} (G \cap \rand B_r) < \frac{\delta}{24Cr}$ imply
$$
  \left| \int_{G \cap \rand B_r} g \, \omega^\epsilon \right|
  < \frac \delta 4 \,.
$$ 
Now, since $d\omega^\epsilon = 0$, we have
$d(g  \omega^\epsilon) = dg \wedge \omega^\epsilon$. 
Since the boundary of $G \cap B_r$ is a smooth hypersurface except for
a set of zero $(n{-}1)$-dimensional Hausdorff measure, we can apply 
Stokes's Theorem (see e.g. \cite[14.3]{Si}) to get
$$
 \int_{G \cap B_r} dg \wedge \omega^\epsilon = 
 \int_{\reg\rand G \cap B_r} g \,\omega^\epsilon +  
 \int_{G \cap \rand B_r} g \,\omega^\epsilon \,.
$$
This implies
$$
 \left| \int_{G \cap B_r} dg \wedge \omega^\epsilon \right| 
 \le \left| \int_{\reg\rand G \cap B_r} g \,\omega^\epsilon \right| + 
  \left| \int_{G \cap \rand B_r} g  \,\omega^\epsilon \right|
 <  
   \frac \delta 2
$$
and
\begin{multline*}
  \left| \int\limits_G dg \wedge \omega \right| 
   \le
   \left| \int\limits_{G \cap B_r} dg \wedge \omega^\epsilon \right| + 
   \left|\int\limits_{G \cap B_r} dg \wedge (\omega^\epsilon - \omega) \right| + 
   \left|\int\limits_{G \setminus B_r} dg \wedge \omega\right| 
   < \delta \,.
\end{multline*}

\paragraph*{Acknowledgment:}
We thank Bruce Kleiner for an interesting discussion on this subject.
We thank the referee for careful reading and thoughtful comments.

Authors' address:\\
  Mathematisches Institut der \\
  Universit\"at Freiburg\\ 
  Eckerstr. 1 \\
  79104 Freiburg \\
  Germany\\
  E-Mail: \\
  franz.auer@math.uni-freiburg.de \\
  bangert@mathematik.uni-freiburg.de

\end{document}